\documentclass[11pt]{article}
 \usepackage{graphicx,amsmath,amsfonts,amssymb,color,mathrsfs}

\usepackage{epsfig}

 \newcommand{\QED}{\hfill \thicklines \framebox(6.6,6.6)[l]{}}

\numberwithin{equation}{section}

 \newtheorem{theorem}{Theorem}[section]
 \newtheorem{lemma}{Lemma}[section]

 \newtheorem{corollary}{Corollary}[section]
 \newtheorem{definition}{Definition}[section]

\newcommand{\eqnb}{\begin{eqnarray*}}
\newcommand{\eqne}{\end{eqnarray*}}

\def\beqlb{\begin{eqnarray}}\def\eeqlb{\end{eqnarray}}
 \def\beqnn{\begin{eqnarray*}}\def\eeqnn{\end{eqnarray*}}

\makeatletter
 \newif\if@borderstar
 \def\bordermatrix{\@ifnextchar*{%
 \@borderstartrue\@bordermatrix@i}{\@borderstarfalse\@bordermatrix@i*}%
 }
 \def\@bordermatrix@i*{\@ifnextchar[{\@bordermatrix@ii}{\@bordermatrix@ii[()]}}
 \def\@bordermatrix@ii[#1]#2{%
 \begingroup
 \m@th\@tempdima8.75\p@\setbox\z@\vbox{%
 \def\cr{\crcr\noalign{\kern 2\p@\global\let\cr\endline }}%
 \ialign {$##$\hfil\kern 2\p@\kern\@tempdima & \thinspace %
 \hfil $##$\hfil && \quad\hfil $##$\hfil\crcr\omit\strut %
 \hfil\crcr\noalign{\kern -\baselineskip}#2\crcr\omit %
 \strut\cr}}%
 \setbox\tw@\vbox{\unvcopy\z@\global\setbox\@ne\lastbox}%
 \setbox\tw@\hbox{\unhbox\@ne\unskip\global\setbox\@ne\lastbox}%
 \setbox\tw@\hbox{%
 $\kern\wd\@ne\kern -\@tempdima\left\@firstoftwo#1%
 \if@borderstar\kern2pt\else\kern -\wd\@ne\fi%
 \global\setbox\@ne\vbox{\box\@ne\if@borderstar\else\kern 2\p@\fi}%
 \vcenter{\if@borderstar\else\kern -\ht\@ne\fi%
 \unvbox\z@\kern-\if@borderstar2\fi\baselineskip}%
 \if@borderstar\kern-2\@tempdima\kern2\p@\else\,\fi\right\@secondoftwo#1 $%
 }\null \;\vbox{\kern\ht\@ne\box\tw@}%
 \endgroup
 }
 \makeatother

\title{\Large \bf GTH Algorithm, Censored Markov Chains, and $RG$-Factorization \\
(dedicated to Dr. Winfried Grassmann)}
\author{
Yiqiang Q. Zhao \\
School of Mathematics and Statistics \\
Carleton University, Ottawa, ON Canada K1S 5B6}
\date{January 2021}

\begin{document}
 \maketitle

\begin{abstract}
In this paper, we provide a review on the GTH algorithm, which is a numerically stable algorithm for computing stationary probabilities of a Markov chain. Mathematically the GTH algorithm is an rearrangement of Gaussian elimination, and therefore they are mathematically equivalent. All components in the GTH algorithm can be interpreted probabilistically based on the censoring concept and each elimination in the GTH algorithm leads to a censored Markov chain.  The $RG$-factorization is a counterpart to the LU-decomposition for Gaussian elimination. The censored Markov chain can also be treated as an extended version of the GTH algorithm for a system consisting of infinitely many linear equations. The censored Markov chain produces a minimal error for approximating the original chain under the $l_1$-norm.
\vspace*{5mm}

\noindent \textbf{Keywords:} GTH method; Gaussian elimination; Markov chains; censored Markov chains; $RG$-factorization; stationary probabilities; numerical stable algorithms
\medskip

\noindent \textbf{Mathematics Subject Classification (2000)} 60J10 $\cdot$ 60J22
\end{abstract}

\section{Introduction}

Stationary probabilities are crucial for stationary behavour of stochastic systems, which can be often modeled as a Markov chain. Explicit expressions for the stationary distribution are available for a small set of problems. In most of applications, simulation or numerical computations are main tools for a solution. Therefore, computational methods for Markov chains are very important. Algorithms are usually presented for a finite-state Markov chain, since computers can only deal with finite-many states. When using finite-state Markov chains to approximate an infinite-state Markov chain, convergence and the approximation error are among fundamental questions to answer.

The GTH algorithm discussed in this paper was proposed by Grassmann, Taksar and Heyman in 1985~\cite{GTH:1985}. Then, it became very attracted to many researchers. A list of  publications by 1993, in which the GTH algorithm was used, can be found in Grassmann~ \cite{Grassmann:1993}, including:
Kohlas~\cite{Kohlas:1986},
Heyman~\cite{Heyman:1987},
Heyman and Reeves~\cite{Heyman-Reeves:1989},
Grassmann and Heyman~\cite{Grassmann-Heyman:1990},
Stewart~\cite{Stewart:1993},
O'Cinneide~\cite{O'Cinneide:1993}.
The list of references closely related to the GTH algorithm,  published after 1993, is so large, and only a small sample is included here:
Stewart~\cite{Stewart:1994},
Dayar and Stewart~\cite{Dayar-Stewart:1996},
Grassmann and Zhao~\cite{Grassmann-Zhao:1997},
Sonin and Thornton~\cite{Sonin-Thornton:2001},
Dayar and Akar~\cite{Dayar-Akar:2005},
Hunter~\cite{Hunter:2018}. Therein from the above publications, many more references can be found.

The GTH algorithm is a numerically stable version of Gaussian elimination, or a rearrangement of Gaussian elimination. Through the rearrangement, subtractions are avoided in the algorithm, which are often the reason causing computational instability. The rearrangement makes the elimination process start with the largest state that controls the error in computations. The GTH algorithm possesses a probabilistic interpretation. This interpretation becomes very clear in terms of the censoring concept of the Markov chain. Two important measures, the $RG$-measures, for Markov chains are invariant under the censoring Based on this property one can prove that the GTH algorithm is equivalent to the $RG$-factorization, which is a counterpart to the LU-decomposition for Gaussian elimination.
The censored Markoc chain also serves as a tool to deal with a countable-state Markov chain. The convergence property given in Section~\ref{sec:5} can be used in many cases to construct a modified countable-state Markov chain such that the GTH algorithm can be performed in the approximation. This convergence property, together with the fact that the censored Markov chain provides an approximation with a minimal error under the $l_1$-norm, leads to efficient stable computations using the GTH algorithm.

The rest of the paper is organized as follows: the GTH algorithm is introduced and discussed in Section~\ref{sec:2};
the censored Markov chain is reviewed in Section~\ref{sec:3}, in which we show that each elimination through using the GTH algorithm results in a censored Markov chain;
the $RG$-factorization is presented in Section~\ref{sec:4}, which is a probabilistic counterpart to the LU-decomposition for Gaussian elimination;
the $RG$-factorization can be considered as an extended version of the GTH algorithm for a system consisting of infinitely many linear equations.

\section{GTH algorithm} \label{sec:2}

The GTH algorithm is a numerical algorithm for computing the stationary distribution of a finite-state Markov chain.
Mathematically, GTH algorithm is a rearrangement of Gaussian elimination. The GTH algorithm is numerically stable, since it starts with smallest entities and subtractions are avoided after the rearrangement, while Gaussian elimination can become numerically unstable if the number of states becomes large.  The GTH algorithm also possesses a probabilistic interpretation in terms of the censoring process.

We start the GTH algorithm with introducing the Markov chain.
\begin{definition}[Markov chain]
A discrete time stochastic process $\{ X_n, n = 0, 1, 2, \ldots \}$, where $X_n$ takes values on a finite or countable set, say $S=\{1, 2,
\ldots \}$, referred to as the state space,  is called a Markov chain if the following Markovian property holds:
\[
    P( X_{n+1}=j | X_n =i, X_{n-1}=i_{n-1}, \ldots, X_1 = i_1, X_0=i_0 )= p_{i,j}
\]
for all states $i_0$, $i_1$, \ldots, $i_{n-1}$, $i$, $j$ and all $n \geq 0$. The Markov chain is called a finite-state Markov chain if $S$ is finite.
\end{definition}

We consider a finite-state Markov chain in this section with the state space
\[
    S=\{1, 2, \ldots, N\}
\]
and the probability transition matrix
\[
    P = \left [  \begin{array}{cccc}
p_{1,1} & p_{1,2} & \ldots \ldots & p_{1,N} \\
p_{2,1} & p_{2,2} & \ldots \ldots & p_{2,N} \\
   \vdots & \vdots & \vdots  & \vdots \\
p_{N,1} & p_{N,2} & \ldots \ldots & p_{N,N} \end{array}
       \right ].
\]
For a finite-state Markov chain, it is well-known that if $P$, or the Markov chain, is irreducible then there exists a unique stationary probability vector (distribution)
\[
    \pi=(\pi_1, \pi_2, \ldots, \pi_N)
\]
satisfying
\[
    \pi = \pi P \quad \text{with} \quad \sum_{n=1}^N \pi_n =1.
\]
Writing out in detail, the equation $\pi = \pi P$ is equivalent to the following $N$ equations (referred to as steady-state, or equilibrium, or stationary equations):
\begin{equation} \label{eqn:se}
    \pi_j = \sum_{i=1}^N \pi_i p_{i,j}, \quad j =1, 2, \ldots, N.
\end{equation}
Since all probabilities $\pi_n$ sum up to one, or $\sum_{n=1}^N \pi_n =1$, one of the above $N$ equations is redundant.

The focus of the GTH algorithm is to numerically compute the stationary probability vector. The basic GTH algorithm consists of two portions: forward eliminations and back substitutions.

\textbf{Forward eliminations:}  First, use the last equation from the system of stationary equations in (\ref{eqn:se}) to eliminate $\pi_N$ to have a system of equations for $N-1$ unknowns with the new coefficients denoted by $p_{i,j}^{N-1}$, where $i, j = 1, 2, \ldots, N-1$, given by
\[
    p_{i,j}^{N-1} = p_{i,j}^N + \frac{p_{i,N}^N p_{N,j}^N}{\sum_{k=1}^{N-1} p_{N,k}^N},
\]
where $p_{i,j}^N = p_{i,j}$. We then have
\[
    \pi_j = \sum_{i=1}^{N-1} \pi_i p_{i,j}^{N-1}, \quad j =1, 2, \ldots, N-1.
\]

It can be directly verified that
\begin{equation} \label{eqn:N-1}
    P^{N-1} \stackrel{\triangle}{=}  \left [  \begin{array}{cccc}
p_{1,1}^{N-1} & p_{1,2}^{N-1} & \ldots \ldots & p_{1,N-1}^{N-1} \\
p_{2,1}^{N-1} & p_{2,2}^{N-1} & \ldots \ldots & p_{2,N-1}^{N-1} \\
   \vdots & \vdots & \vdots& \vdots \\
p_{N-1,1}^{N-1} & p_{N-1,2}^{N-1} & \ldots \ldots & p_{N-1,N-1}^{N-1} \end{array}
       \right ]
\end{equation}
is a stochastic matrix, or a new Markov chain.

We then repeat the above elimination process to eliminate $\pi_{N-1}$, $\pi_{N-2}$, \ldots, $\pi_{n}$ to have the following coefficients
\[
    p_{i,j}^{n-1} = p_{i,j}^n + \frac{p_{i,n}^n p_{n,j}^n}{\sum_{k=1}^{n-1} p_{n,k}^n}
\]
and $\pi_j$ satisfy the following equation:
\begin{equation} \label{eqn:se-j}
    \pi_j = \sum_{i=1}^{n-1} \pi_i p_{i,j}^{n-1}, \quad j =1, 2, \ldots, n-1.
\end{equation}
The matrix of these coefficients:
\begin{equation} \label{ean:n-1}
    P^{n-1} \stackrel{\triangle}{=}  \left [  \begin{array}{cccc}
p_{1,1}^{n-1} & p_{1,2}^{n-1} & \ldots \ldots & p_{1,n-1}^{n-1} \\
p_{2,1}^{n-1} & p_{2,2}^{n-1} & \ldots \ldots & p_{2,n-1}^{n-1} \\
   \vdots & \vdots & \vdots& \vdots \\
p_{n-1,1}^{n-1} & p_{n-1,2}^{n-1} & \ldots \ldots & p_{n-1,n-1}^{n-1} \end{array}
       \right ]
\end{equation}
defines a Markov chain.

After $\pi_2$ is eliminated, we reach a Markov chain of size 1, or $p_{1,1}^1=1$, which means that we have
\[
    \pi_1 = \pi_1 \cdot 1.
\]

It is worthwhile to mention that the forward elimination starts from $\pi_N$, then $\pi_{N-1}$, and move to $\pi_n$ with smaller values of $n$. This is important for the control of the computation errors, since for large $n$, tail probabilities are smaller. Also, $1-p_{n,n}^n$ is replaced by $\sum_{k=1}^{n-1} p_{n,k}^n$ to avoid the subtraction.

\textbf{Back substitutions:}
To find the solution for $\pi_j$, the GTH performs the back substitution. Define
\[
    r_1=1, \quad r_j = \pi_j/\pi_1, \quad j=2, 3, \ldots, N.
\]
Then, $r=(r_1, r_2, \dots, r_N)$ satisfy equation (\ref{eqn:se}) and also (\ref{eqn:se-j}). First, substitute $r_1$ in the system with coefficients $p_{i,j}^2$ to have
\begin{align*}
    r_1 =& r_1 p_{1,1}^2 + r_2 p_{2,1}^2, \\
     r_2 = & r_1 p_{1,2}^2 + r_2 p_{2,2}^2.
\end{align*}
Since one of the two equations is redundant, we take the second one to have
\[
    r_2 = r_1 \frac{p_{1,2}^{2}}{p_{2,1}^2},
\]
and use $1-p_{2,2}^2 = p_{2,1}^2$ to avoid the substraction in the algorithm.
Repeat the above back substitution process until $\pi_{j-1}$ has been substituted to have
\[
    r_j = \sum_{i=1}^{j-1} r_i \frac{p_{i,j}^{j}}{\sum_{k=1}^{j-1} p_{j,k}^j}, \quad j = 1, 2, \ldots, N.
\]

Recall that $r=(r_1, r_2, \dots, r_N)$ and $\pi=(\pi_1, \pi_2, \ldots, \pi_N)$ are different only by a constant, and $\pi$ is a probability vector, we can easily normalize $r$ to have
\begin{equation} \label{eqn:rj}
    \pi_j = \frac{r_j}{\sum_{k=1}^N r_k}, \quad j = 1, 2, \ldots, N.
\end{equation}

It is not difficult to see the mathematical equivalence between the GTH algorithm and Gaussian elimination. As indicated earlier, Gaussian elimination is usually numerically unstable when $N$ is large, say $10,000$ or larger, while the GTH algorithm is very stable.

\section{Censored Markov chains} \label{sec:3}

The GTH algorithm has probabilistic interpretations. During the forward elimination, each step results in a new Markov chain with the state space one state fewer than the previous state space. In fact, each of these Markov chains is a so-called censored Markov chain, which will be discussed in this section.

\begin{definition}[Censored process]
Let $\{X_n: n =0, 1, 2, \ldots \}$ be a countable-state Markov chain with state space
$S=\{1, 2, \ldots \ldots \}$ and probability transition matrix $P$.
Let $E$ be a non-empty subset of $S$. Suppose that $n_k$ is the $k$th time at which the original process $X_n$ visits the subset $E$. Then, the censored stochastic process $\{X_k^E: k =0, 1, 2, \ldots \}$ is defined by $X_k^E=X_{n_k}$, or the value of the new process $X_k^E$ at time $k$ is equal to the value of the original process $X_n$ at its $k$th time of visiting the subset $E$.
\end{definition}

The censored process is also referred to as a watched process since it is obtained by watching $X_n$ only when it is in $E$. It is also referred to as an embedded Markov chain since the time (or the state space) of the censored process is embedded in the time (or the state space) of the original process.
The following lemma is a summary of some basic properties of the censored process.
\begin{lemma} \label{lem:censored}
\textbf{(i)} The censored process $X_k^E$ is also a Markov chain. If the probability transition matrix $P$ of the Markov chain $X_n$ is partitioned according to $E$ and its complement $E^c$:
\[
    P = \bordermatrix[{[]}]{%
    & E & E^c \cr
E & T & U \cr
E^c & D & Q
},
\]
then the probability transition matrix of the censored Markov chain $X_k^E$ is given by
\[
    P^E = T + U \hat{Q} D,
\]
where $\hat{Q} = \sum_{n=0}^{\infty} Q^n$ is the minimal inverse of $I-Q$.

\textbf{(ii)}
The Markov chain $P$ is irreducible, if and only if for every subset $E$ of the state space, the censored Markov chain $P^E$ is irreducible. If $P$ has a unique stationary probability vector $\pi=(\pi_1, \pi_2, \ldots)$, then the stationary probability vector $\pi^E = \big (\pi_j^E \big )_{j \in E}$ of the censored Markov chain is given by
\[
	\pi^E_j = \frac{\pi_j}{\sum_{k \in E} \pi_k},\quad j \in E.
\]

\textbf{(iii)} If $E_1$ and $E_2$ are two non-empty subsets of the state space $S$ and $E_2$ is a subset of $E_1$, then
\[
    P^{E_2} = (P^{E_1})^{E_2}.
\]
\end{lemma}

The concept of the censored Markov chain was first introduced and studied by L\'{e}vy~\cite{Lev:1951,Lev:1952,Lev:1958}. It was then used by Kemeny, Snell and Knapp~\cite{Kem66}
for proving the uniqueness of the invariant vector for a recurrent countable-state Markov chain. This embedded Markov chain was an approximation tool in the book by Freedman~\cite{Fre:1983} for countable-state Markov chains. When the censored Markov chain is used to approximate the stationary distribution, Zhao and Liu~\cite{Zha:96} proved that it has the smallest error in $l_1$-norm among all possible approximations.

Now, we discuss the connection between the GTH algorithm and the censored Markov chain. First, it is easy to check that if we let $E_n=\{1, 2, \ldots, n \}$, then the Markov chain
$P^{N-1}$ in (\ref{eqn:N-1}) is the censored Markov chain $P^{E_{N-1}}$, and $P^n$ in (\ref{ean:n-1}) is the censored Markov chain  $P^{E_{n-1}}$ according to Lemma~\ref{lem:censored}. The expression for the stationary distribution given in (\ref{eqn:rj}) is an immediate consequence of Lemma~\ref{lem:censored}-(ii). More probabilistic interpretations for the GTH algorithms can be provided:
\begin{description}
\item[(1)] $(1-p_{n,n}^{n})^{-1} = \big (\sum_{k=1}^{n-1} p_{n,k}^n \big )^{-1}$ is the expected number of visits to state $n$ before entering the censored set $E_{n-1}$ given that the process started in $n$.
\item[(2)]
\[
    \frac{p_{i,n}^n }{\sum_{k=1}^{n-1} p_{n,k}^n}, \quad i < n,
\]
is the expected number of visits to state $n$ before returning to $E_{n-1}$ given that the process started in
    state $i < n$.
\item[(3)]
\[
    \frac{p_{n,j}^n}{\sum_{k=1}^{n-1} p_{n,k}^n}, \quad j < n,
\]
is the probability that upon entering $E_{n-1}$
    the first state visited is $j < n$, given that the process
    started in state $n$.
\item[(4)]
\[
    \frac{p_{i,n}^n p_{n,j}^n}{\sum_{k=1}^{n-1} p_{n,k}^n}, \quad i, j < n,
\]
is the probability that upon returning to $E_{n-1}$
    the first state visited is $j$, given that the process
    started in state $i < n$.
\item[(5)]
\[
    p_{i,j}^{n-1} = p_{i,j}^n + \frac{p_{i,n}^n p_{n,j}^n}{\sum_{k=1}^{n-1} p_{n,k}^n}
\]
is the transition probability from $i$ to $j$ of the censored Markov chain $P^{E_{n-1}}$.
\end{description}

\section{$RG$-factorization} \label{sec:4}

Mathematically, the GTH algorithm is equivalent to Gaussian elimination, which in turn is equivalent to an LU-factorization (or UL-factorization). In this section, we discussion the $RG$-factorization for Markov chains and show that the GTH algorithm is equivalent to the $RG$-factorization.

The $RG$-factorization discussed here is one of the versions of the so-called Wiener-Hopf-type factorization. This version of factorization is given in terms of the dual measures, the $RG$-measures, of the Markov chain. People who are interested in this topic could refer to the literature references, including
Heyman~\cite{Heyman:1995},
Zhao, Li and Braun~\cite{ZLB:1997},
Zhao~\cite{Zhao:2000},
Li and Zhao~\cite{Li-Zhao:2002,Li-Zhao:2002b,Li-Zhao:2004}.

Consider an irreducible countable-state Markov chain with its probability transition matrix $P$ on the state space
\[
    S=\{1, 2, 3, \ldots \},
\]
given by
\begin{equation} \label{eqn:Pij}
    P =\left [
\begin{array}{cccc}
p_{1,1} & p_{1,2} & p_{1,3} & \cdots \\
p_{2,1} & p_{2,2} & p_{2,3} & \cdots \\
p_{3,1} & p_{3,2} & p_{3,3} & \cdots \\
\vdots & \vdots & \vdots & \ddots%
\end{array}
\right ] .
\end{equation}

Define a pair of dual measures as follows:
    for $1 \leq i \leq j$, define $r_{i,j}$ to be the expected number of visits to state $j$ before hitting any state $< j$, given that
the process starts in state $i$; and
    for $i > j \geq 1$, define $g_{i,j}$ to be the probability of hitting state $j$ for the first time, given that the process starts in state $i$.

One of the most important properties for the $RG$-measures is the invariance under censoring, which is stated in the following theorem.
\begin{theorem}[Invariance of $RG$-measures, Theorem~4 in \cite{ZLB:2003}] \label{theorem1}
    For the Markov chain given in (\ref{eqn:Pij}) with its $RG$-measures $r_{i,j}$ and $g_{i,j}$. Let $r_{i,j}^n$ and $g_{i,j}^n$ be the $RG$-measures for the censored Markov chain $P^{E_n}$. Then, for given $1=i<j$, or $1 \leq i \leq j$,
\begin{equation} \label{Eqn:13}
    r^n_{i,j} = r_{i,j}, \quad \text{for all $n \geq j$},
\end{equation}
and for given $1 \leq j < i$,
\begin{equation}
    g^n_{i,j} = g_{i,j}, \quad \text{for all $n \geq i$.}
\end{equation}
\end{theorem}

The $RG$-factorization of the Markov chain $P$ is given in the following theorem.
\begin{theorem}[$RG$-factorization, Theorem~13 in \cite{Zhao:2000}] \label{The:RGfactor}
    For the Markov chain defined by (\ref{eqn:Pij}), we have
\begin{equation}
    I-P = \left[ I-R_{U} \right] \left [ I- \Psi_{D}\right] \left[ I-G_{L} \right],  \label{gRGdecomp}
\end{equation}
where
\begin{equation*}
    R_{U} =\left [
\begin{array}{ccccc}
0 & r_{1,2} & r_{1,3} & r_{1,4} & \cdots \\
& 0 & r_{2,3} & r_{2,4} & \cdots \\
&  & 0 & r_{3,4} & \cdots \\
&  &  & 0 & \cdots \\
&  &  &  & \ddots
\end{array}
\right ],
\end{equation*}
\begin{equation*}
    \Psi_{D} =\text{diag}\left( \psi_{1},\psi_{2},\cdots \right)
\end{equation*}
with $\psi_n=p^n_{n,n}$, and
\begin{equation*}
    G_{L} =\left [
\begin{array}{ccccc}
0 &  &  &  &  \\
g_{2,1} & 0 &  &  &  \\
g_{3,1} & g_{3,2} & 0 &  &  \\
g_{4,1} & g_{4,2} & g_{4,3} & 0
&  \\
\vdots & \vdots & \vdots & \vdots & \ddots%
\end{array}
\right ].
\end{equation*}
\end{theorem}

The $RG$-factorization of the GTH algorithm is an immediate consequence of the above theorem and the invariance of the $RG$-measures under censoring.
\begin{corollary}[$RG$-factorization of GTH algorithm]
The GTH algorithm is equivalent to the following $RG$-factorization:
\begin{equation*}
    I-P = \left[ I-R_{U} \right] \left [ I- \Psi_{D}\right] \left[ I-G_{L} \right],
\end{equation*}
where
\begin{equation*}
    R_{U} =\left [
\begin{array}{cccccc}
0 & r_{1,2} & r_{1,3} & r_{1,4} & \cdots & r_{1,N}\\
& 0 & r_{2,3} & r_{2,4} & \cdots & r_{2,N} \\
&  & 0 & r_{3,4} & \cdots & r_{3,N} \\
&  &  & 0 & \cdots & \vdots \\
&  &  &  & \ddots & r_{N-1,N} \\
&  &  &  & & 0
\end{array}
\right ],
\end{equation*}
\begin{equation*}
    \Psi_{D} =\text{diag}\left( \psi_{1},\psi_{2},\cdots, \psi_N \right),
\end{equation*}
and
\begin{equation*}
    G_{L} =\left [
\begin{array}{cccccc}
0 &  &  &  &  \\
g_{2,1} & 0 &  &  &  \\
g_{3,1} & g_{3,2} & 0 &  &  \\
g_{4,1} & g_{4,2} & g_{4,3} & 0 \\
\vdots  &  \vdots & \vdots  & \vdots & \ddots \\
g_{N,1} & g_{N,2} & g_{N,3} & \cdots & g_{N,N-1}& 0 \\
\end{array}
\right ].
\end{equation*}
Specifically, for $n =1, 2, \ldots, N$,
\begin{align*}
    r_{i,n} = & \frac{p_{i,n}^n }{\sum_{k=1}^{n-1} p_{n,k}^n}, \quad i = 1,, 2, \ldots, n-1, \\
    g_{n,j} = &     \frac{p_{n,j}^n}{\sum_{k=1}^{n-1} p_{n,k}^n}, \quad j = 1, 2, \ldots, n-1, \\
    \psi_n= & p^n_{n,n}.
\end{align*}
\end{corollary}

\section{Extending GTH to countable-state Markov chains} \label{sec:5}

Gaussian elimination is a method in linear algebra for solving a system of finitely many linear equations. Since linear algebra deals with linear spaces with a finite dimension, there is no Gaussian elimination version for a system consisting of infinitely many linear equations. The $RG$-factorization, together with the censored Markov chain, can be treated as an extended version of the GTH algorithm, since the UL-factorization formally allows us to perform forward eliminations and back substitutions to compute the stationary vector $\pi$. However, in order to practically start the elimination, we need a start state, not the infinite, which leads to various truncation methods.

Assume that for a large $N$, we partition the transition matrix $P$ according to $E_N=\{1, 2, \ldots, N\}$ and $E^c_N$:
\begin{equation} \label{eqn:partition}
    P =\left [
\begin{array}{cccc}
p_{1,1} & p_{1,2} & p_{1,3} & \cdots \\
p_{2,1} & p_{2,2} & p_{2,3} & \cdots \\
p_{3,1} & p_{3,2} & p_{3,3} & \cdots \\
\vdots & \vdots & \vdots & \ddots%
\end{array}
\right ]  = \bordermatrix[{[]}]{%
    & E_N & E_N^c \cr
E_N & T_N & U \cr
E_N^c & D & Q
}
\end{equation}
An augmentation is a method using a non-negative matrix $A_N$ such that
\begin{equation} \label{eqn:augmentation}
    \tilde{P}_N = T_N + A_N
\end{equation}
is stochastic. Popular augmentations include the censored Markov chain, the last column augmentation (add the missing probabilities to the last column), the first column augmentation (add the missing probability to the first column), and more generally (than the last and first column augmentations), the linear augmentation (add the missing probabilities linearly to the first $N$ columns).

Assume that $\tilde{P}_N$ is irreducible. Under this condition, let $\tilde{\pi}^N$ be the unique stationary probability vector of the augmented Markov chain. We are interested in the convergence of $\tilde{\pi}^N$ to $\pi$ as $N \to \infty$, and the error between $\tilde{\pi}^N$ and $\pi$ when the convergence is established.
These two questions are fundamental in the area of approximating a countable-state Markov chain by finite-state Markov chains. The book of Freedman~\cite{Fre:1983} was devoted to the first question using censored Markov chains. Based on \cite{Fre:1983}, we can develop the following two convergence properties.
\begin{lemma}[Lemma~3 in \cite{Zhao:2000}] \label{lemma3}
Let $P$ be a transition matrix with state space $S$
and let $E_n$ for $n = 1, 2, \ldots$ be a sequence of subsets of $S$
such that $E_n \subseteq E_{n+1}$ and $lim_{n \to \infty} E_n = S$.
Then, for any $i, j \in E_n$,
$\lim_{n \to \infty} P_{i,j}^{E_n} = P_{i,j}$.
\end{lemma}
This lemma says that if the sequence $E_n$ (not necessarily equal to $\{1, 2, \ldots, n\}$) of subsets of the state space converges to the state space $S$, then the sequence of the censored Markov chains converges to the original chain.

\begin{theorem}[Theorem~7 in \cite{Zhao:2000}] \label{the:7}
    Let $P=(p_{i,j})_{i,j = 1, 2, \ldots}$ be the transition matrix of a recurrent
Markov chain on the positive integers. For an integer $\omega > 0$,
let  $P(\omega )= (p_{i,j}(\omega ))_{i,j =1, 2, \ldots}$ be a matrix
such that
\[
    p_{i,j}(\omega ) = p_{i,j}, \;\;\;  \mbox{for } i,j \leq \omega
\]
and $P(\omega )$ is either stochastic or substochastic matrix.
For any fixed $n \geq 0$, let $E_n = \{ 1, 2,  \ldots, n \}$ be
the censoring set. Then,
\begin{equation}
    \lim_{\omega \to \infty} P^{E_n}(\omega ) = P^{E_n}.
\end{equation}
\end{theorem}

This theorem provides us with many options to construct an infinite stochastic matrix $P(\omega)$ with the same northwest corner $T_\omega$ such that the censored Markov chain $P^{E_N}(\omega)$ can be easily obtained. Then, we can apply the GTH algorithm to the finite-state Markov chain $P^{E_N}(\omega)$ to compute the stationary probability vector $\pi^N$ for the censored Markov chain $P^{E_N}(\omega)$. According to the above theorem (Theorem~\ref{the:7}), $\pi^N$ is an approximation to $\pi$. This procedure also results in an approximation with an ``approximate'' minimal error in the sense of $l_1$-form based on the main result in Zhao~\cite{Zha:96}. We provide a brief discussion here.

Consider an irreducible augmentation $\tilde{P}_N$ defined in (\ref{eqn:augmentation}) with the unique stationary probability vector $\tilde{\pi}^N$. Suppose that the augmentation is convergent in the sense that
\[
    \lim_{N \to \infty}  \tilde{\pi}_j^N = \pi_j, \quad j =1, 2, \ldots.
\]
We are interested in the error between $\tilde{\pi}^N$ and $\pi$. The $l_1$-norm of this error is defined by
\[
	l_1(N, \infty) = \sum_{k=1}^{N} \big | \tilde{\pi}^N_k - \pi_k \big |+ \sum_{k=N+1}^{\infty} \pi_k.
\]
It is worthwhile to comment here that not all augmentations are convergent. For example, the popular last column augmentation may not be convergent (see for example, \cite{Gib:1987a}).

The following result guarantees a minimal error sum for the censored Markov chain.
\begin{theorem}[Best augmentation, \cite{Zha:96}]
    The censored Markov chain is an augmentation method such that the error sum $l_1(K, \infty)$ is the minimum.
\end{theorem}

The first column augmentation is the worst under the $l_1$-norm and the last column augmentation, if it is convergent, is the best under the $l_\infty$-norm (see also \cite{Li-Zhao:2000}). Other references on augmentations include \cite{S80,T98,M17,ll18}, and references therein.

\section{Concluding words} \label{sec:6}

This review paper is dedicated to Dr.~Winfried Grassmann, who is my Ph.D. supervisor, who directed me to the area of queueing theory and applied/computational probability. The GTH algorithm is one of his celebrated contributions to applied probability, and it is now a standard textbook content for computations of Markov chains.

\vspace*{3mm}
\noindent \textbf{Acknowledgements:} The author acknowledges that this work was supported in part through a Discovery Grant of NSERC.


\end{document}